\documentclass[11pt]{article}
\usepackage{amssymb,sans}
\usepackage{graphicx,subfigure}

\begin{document}

\newcommand{\ea}{et al.}
\newcommand{\be}{\begin{equation}}
\newcommand{\ee}{\end{equation}}
\baselineskip 17pt
\parindent 10pt
\parskip 8pt

\centerline{\Large\bf Of Bombs and Boats and Mice and Men}

\centerline{\large\bf A random tour through some scaling laws}

\vskip 0.1in
\centerline{\bf Niall MacKay}

\noindent{\bf I.} I'll begin with a true  story which has become semi-mythical.  Shortly after World War Two, the US Atomic Energy Commission released a film of the 1945 Trinity atomic bomb test. The energy yield remained secret, having been estimated only with some difficulty. So the Americans were most surprised when the British fluid dynamicist Geoffrey Taylor published, in 1950, an accurate estimate merely by studying the AEC's pictures \cite{Taylor}. Embellishments of the tale have the Americans wondering how on Earth he did it, with the CIA visiting his house in the middle of the night to search through his papers.

\vskip 0.2in
\centerline{\includegraphics[keepaspectratio=true,width=3in]{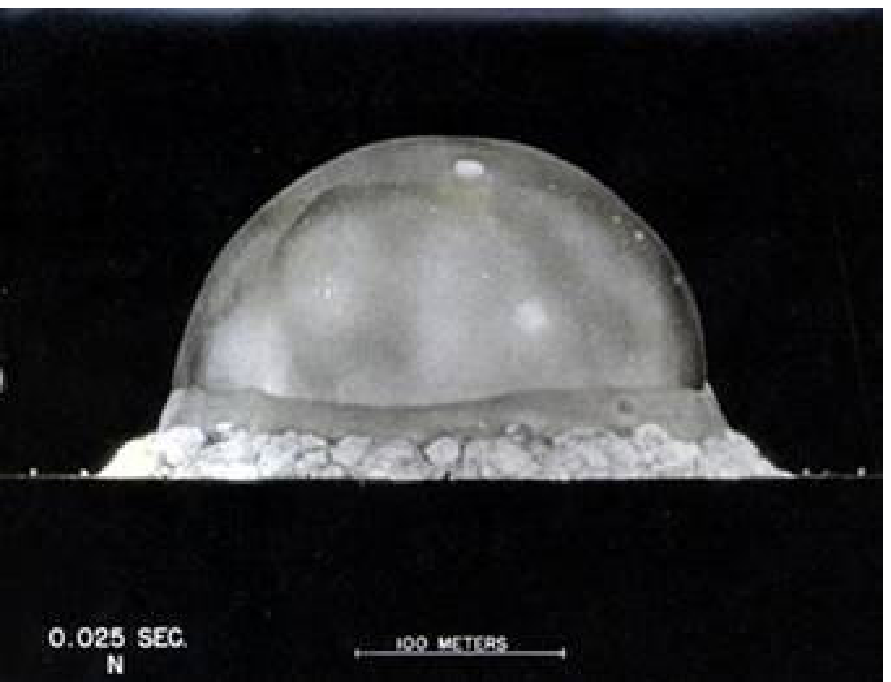}}
\centerline{Figure 1: the Trinity atomic bomb test.}

The way the story is often told nowadays, his technique was simple dimensional analysis. Suppose that, at time $t$, the radius $r$ of the blast wave depends only on $t$ (the time elapsed since the explosion), on the air density $\rho$, and on the energy $E$ released in the blast. (That is, once the blast wave is propagating it knows nothing of the nature of the explosion which caused it.) Energy is a mass multiplied by the square of a speed, which we write $[E]=ML^2T^{-2}$, while for the density $[\rho]=ML^{-3}$. So to cancel the masses, we must combine $E$ and $\rho$ as $E/\rho$. This has dimensions $[E/\rho]=L^5T^{-2}$, so that to form a length from this and $t$ we must take $(Et^2/\rho)^{1/5}$, and finally
\begin{equation}\label{bomb}
r=C\left( {Et^2\over \rho} \right)^{1\over 5}\,,
\end{equation}
where $C$ is a dimensionless constant, universal for such waves, which can be estimated from conventional explosive blasts. Used with pictures such as that in Fig.\ 1, which helpfully gives length and time scales, this formula enabled Taylor to estimate $E$.

This calculation is striking especially for the unexpected appearance of the one-fifth power---after all, when we clap our hands we are used  to the effects propagating very differently, with a fixed `speed of sound' $v$, and $r\propto vt$. It is often used to illustrate the power of dimensional analysis, but this is really another myth, potentially more pernicious than the more exciting embellishments, for the dimensional result only emerged as part of a deeper analysis which included both the mathematics of self-similar solutions of the relevant partial differential equations (PDEs) and a careful treatment of the physics involved. The dimensional derivation is true, but hardly independent.\footnote{In fact Taylor had been working on the mathematics of blast waves throughout the war, publishing his results in a special, classified volume of the {\em Proceedings of the Royal Society}. Thanks to Martin Smith for drawing my attention to this volume.}

Another classic example is the problem of roasting times. How much longer does it take to roast a 5kg turkey than a 1kg pheasant? Most recipes give a linear formula: so many minutes per unit of mass, plus a fixed time. However, this is only an approximation to a scaling which can be discovered by dimensional analysis. Let us assume that the birds are similar in shape, roasted at the same temperature, with similar flesh, and that we wish their centres to reach the same temperature to be properly cooked. Then the cooking time $t$ depends only on the typical distance $l$ from surface to centre and on one physical parameter, the thermal diffusivity $\kappa$. This is associated with the one microscopic physical law which we need to solve the problem, the `heat equation': that the rate of change of temperature is proportional (via $\kappa$) to its Laplacian (of second spatial derivatives), so that $[\kappa]=L^2T^{-1}$. (We don't need to consider precisely how $\kappa$ depends on the thermal conductivity, density and specific heat capacity of the flesh.) This immediately gives us
 that $t\sim l^2/\kappa$. Since the bird's mass $m\sim l^3$, we have $t=C' m^{2/3}/\kappa$ for some dimensionless constant $C'$. So, if the roasting time for the pheasant is $t_p\simeq 1$hr, then
$$ {t_{\rm turkey} \over t_{\rm pheasant}} = \left({m_{\rm turkey} \over m_{\rm pheasant}}\right)^{2/3}
$$ and the turkey's roasting time is $5^{2/3}\simeq 3$ hours.\footnote{We read $y\sim x^\beta$ as `$y$ scales as $x^\beta$', and mean (with no more precision than is necessary for our purposes) that $yx^{-\beta}$ is approximately constant. If the reader baulks at the use of $l$ and the assertion that $m\sim l^3$, we note that one could instead use $m$ and density $\rho$, with $[\rho]=ML^{-3}$.}

Dimensional analysis, then, is fundamentally bound up with the concept of {\em scaling}, which is crucial in applying mathematics. For example, suppose we wish to scale up a laboratory chemical experiment to an industrial process. If the reaction produces heat then we see immediately that the apparatus may not simply be self-similarly scaled up, for then the extra heat produced (in proportion to the volume of reagents) will not be matched by the ability to dissipate it safely (in proportion to the surface area of the vessel). Rather we should identify the dimensionless numbers which characterize  the problem and  keep these invariant under the scaling. In the examples above, there is only one dimensionless number, $Et^2/\rho r^5$ for the bomb blast and $\kappa t / l^2$ for the roasting times.  If there were more, we could write down in scale-free form the equations which relate them.

To keep a dimensionless number constant is precisely to impose a power-law scaling, as in (\ref{bomb}). Indeed for scale-invariance $r$ {\em must} scale as a simple power of each of $E$, $\rho$ and $t$, for this is the only functional form which is invariant under changes of scale (or, equivalently, of units). With $r=K t^{2/5}$ (taking $K=C(E/\rho)^{1/5}$ to be constant), multiplying $t$ and $r$ by constant factors merely changes $K$; it does not change the form of the equation, the power law. It was precisely this similarity property of the blast wave that enabled Taylor to simplify the solution of his PDE.

The basic picture we have described holds in all applied mathematics, and dimensional analysis is an essential tool, not only for all natural scientists but for anyone who applies mathematics. Even money, for example, can be considered a primitive dimensionful quantity, with its own units (of currency). Yet the subject is taught only patchily to users of maths: in English schools it appears (as far as I am aware) only in the more advanced A-level mechanics modules. There is a particular danger, as we shall see in Sect.\ IV, in the practice of not including the units in the algebraic variables----that is, for example, in writing a speed as $v$ ms$^{-1}$ with $v$ a (dimensionless) number, rather than simply as $v$, a dimensionful quantity with no intrinsic number but which acquires one when divided by some suitable unit $v_0$, such as $v_0=1$ ms$^{-1}$. It is striking that this is practised in the same examinations which test elementary dimensional analysis---surely a recipe for confusion.

\vskip 0.1in
\noindent{\bf II.}

\centerline{\includegraphics[keepaspectratio=true,
width=120pt]{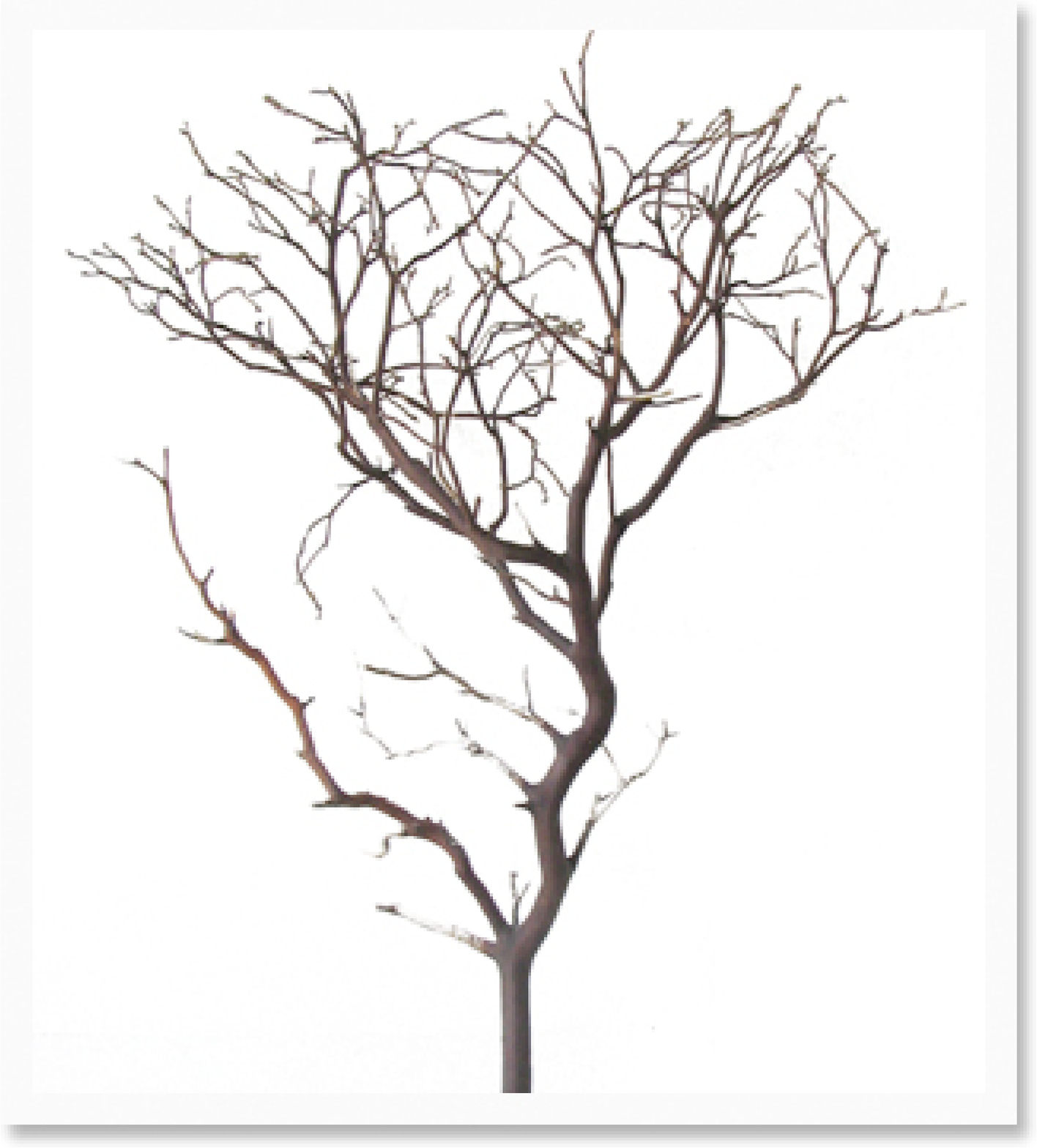}
\hspace*{0.7in}\includegraphics[keepaspectratio=true,
width=117pt]{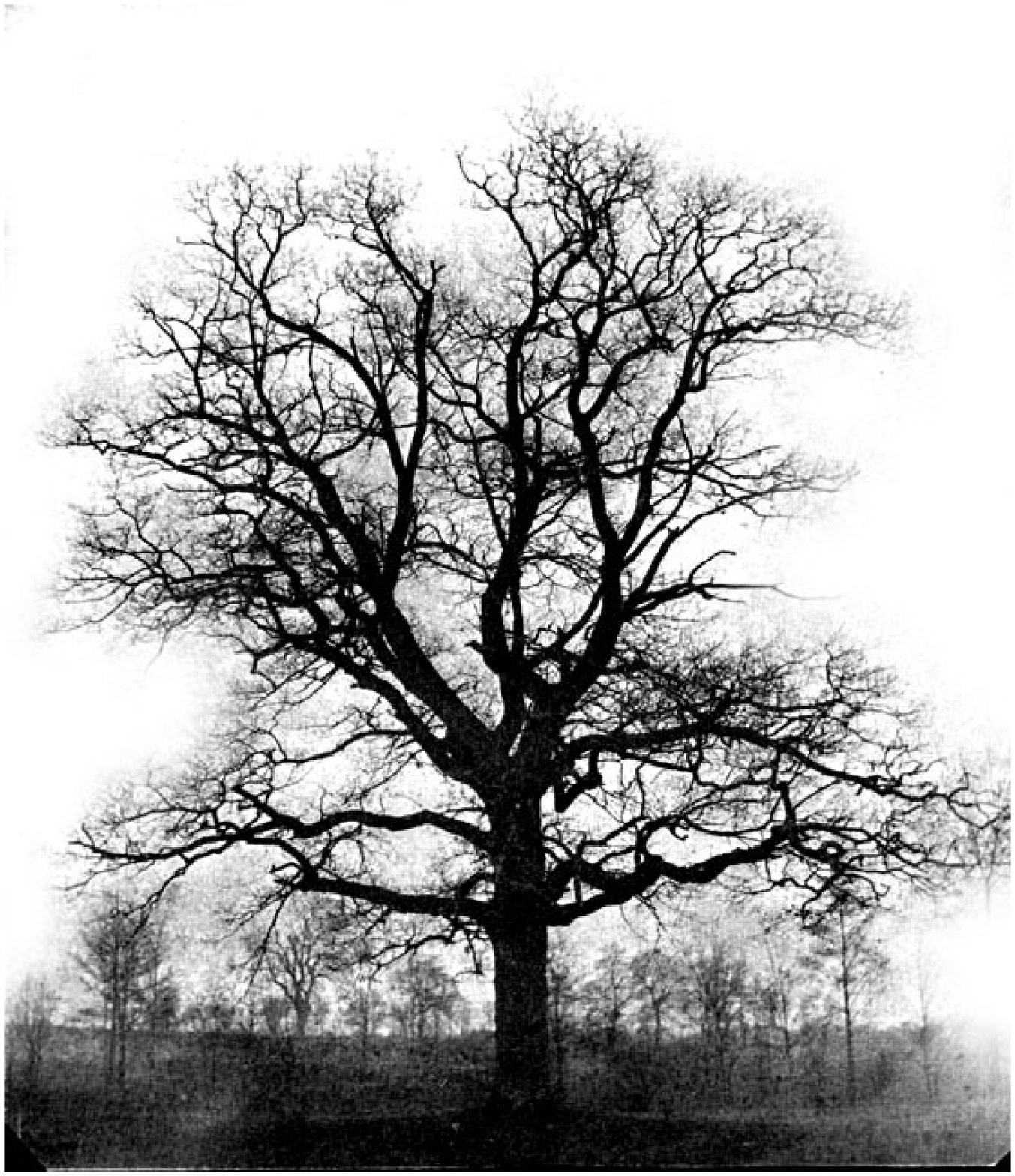}}

\centerline{Figure 2: a 1' twig and a 100' tree.}
\noindent An impression is sometimes given in popular science that scale-invariance is most naturally realized as a property of branching networks---so that a twig, for example, has much the same architecture as a tree (Fig.\ 2). But of course we can find it in any process which is invariant under rescaling of a dimensionful quantity. Restricting for the moment to lengths, how about the pictures in Fig.\ 3? With modern materials and handling gear, even a large yacht can have a single large mainsail, and such `Bermudan sloop' design is approximately invariant over an order of magnitude.

\vskip 0.2in
\centerline{\includegraphics[keepaspectratio=true,
width=120pt]{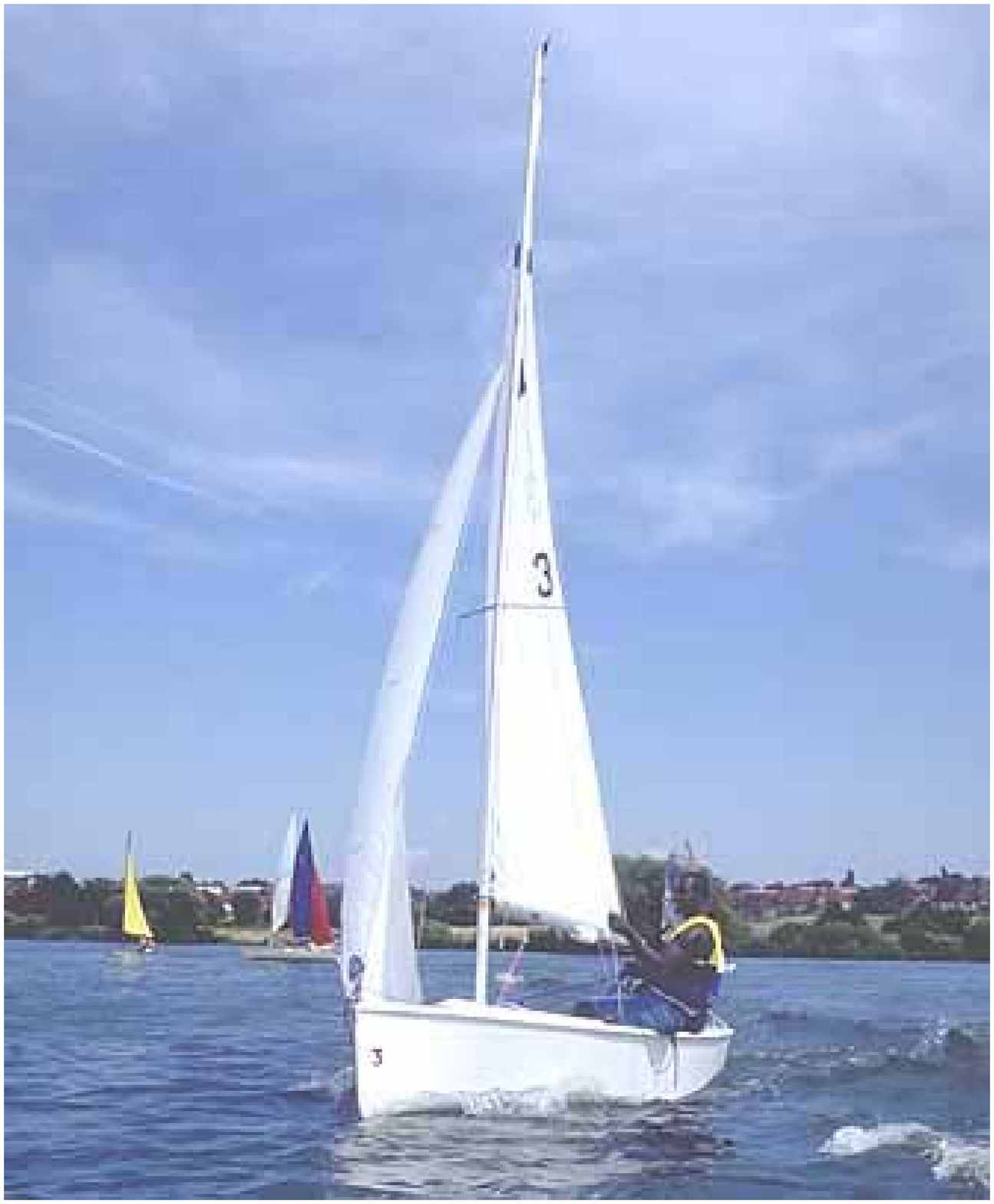}
\hspace*{0.7in}\includegraphics[keepaspectratio=true,
width=117pt]{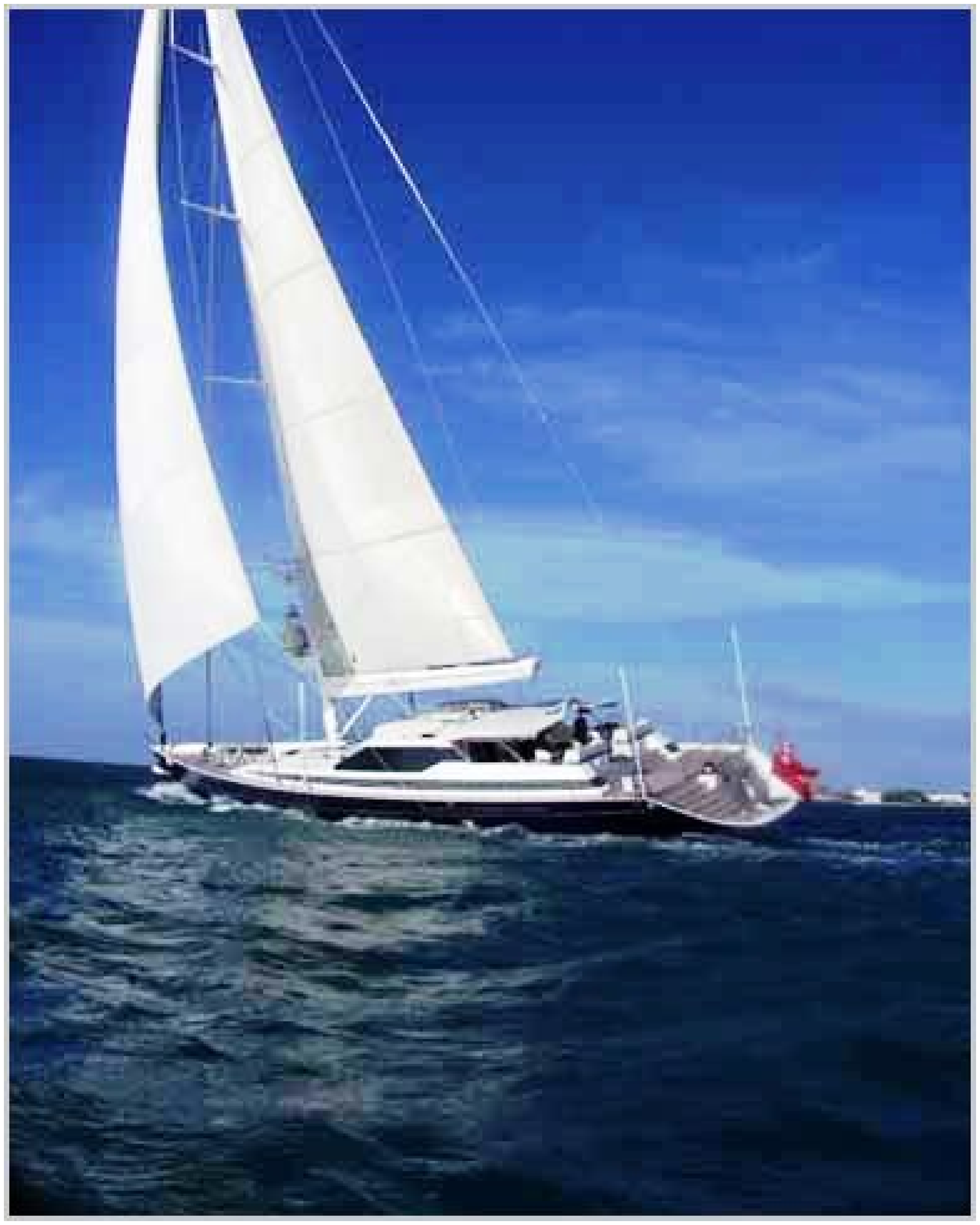}}

\centerline{Figure 3: a 14' dinghy and a 109' yacht.}

So what power laws might we seek for yachts? An obvious one is speed: how does the maximum speed of a yacht scale with its length? Well, {\em if} the speed depends only on the boat's length $l$ and on the acceleration $g$ due to gravity, then the only way to extract a speed from these is as the square root of their product: $v \sim \sqrt{gl}$. This, of course, is where the mystery of the game of dimensional analysis lies: which variables should be invited to the party? Why $g$, but no property of the hull materials or of water? In this case the answer is straightforward: every `displacement' hull, which makes progress by moving through the water, sets up a bow wave, of wavelength proportional to the boat's waterline length, and which it cannot overtake. A surface wave of wavelength $\lambda$ on deep water has speed $\sqrt{g\lambda/2\pi}$, and the result follows.\footnote{For a classic introduction see \cite{Coulson}. The wave speed does not depend on properties of the water, and is typically valid for any incompressible, inviscid fluid. (Note that viscosity would introduce another classic dimensionless quantity, the Reynolds number.) The scaling of ship speed is quite accurate: a 25' yacht has a maximum speed of about 6 knots, and a 600' aircraft carrier, however powerful, has a maximum speed of about 30 knots.} But such questions are generally very subtle \cite{Barenblatt}.

So let's be flippant instead. How fast does a yacht's {\em price} scale with length? There's no underlying physical law here, just the kind of question to which a mathematician at a boat show, astonished by the high prices, wants an answer.

\centerline{\includegraphics[keepaspectratio=true,
width=200pt]{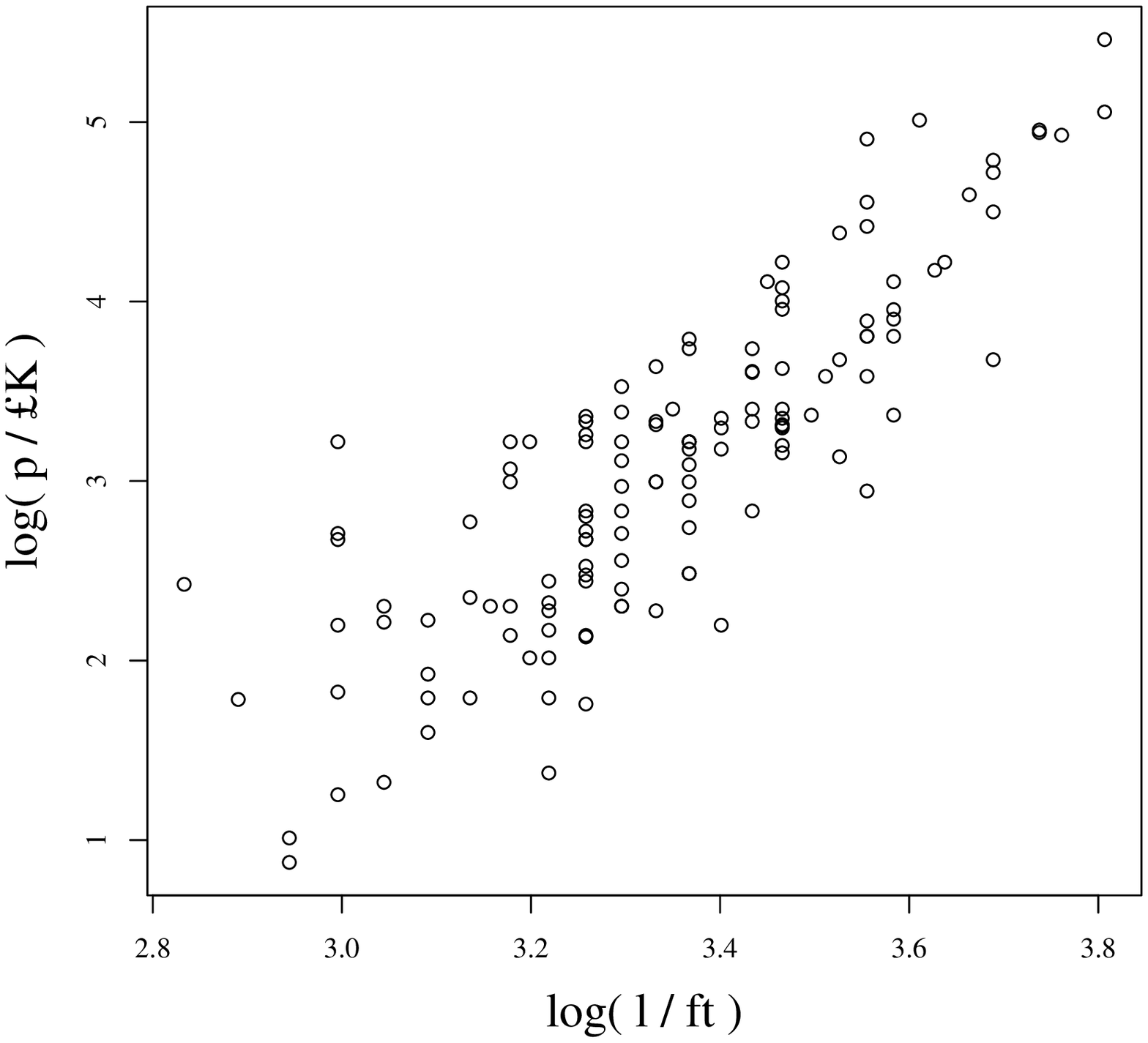}
\hspace{0.6in}\includegraphics[keepaspectratio=true,
width=200pt]{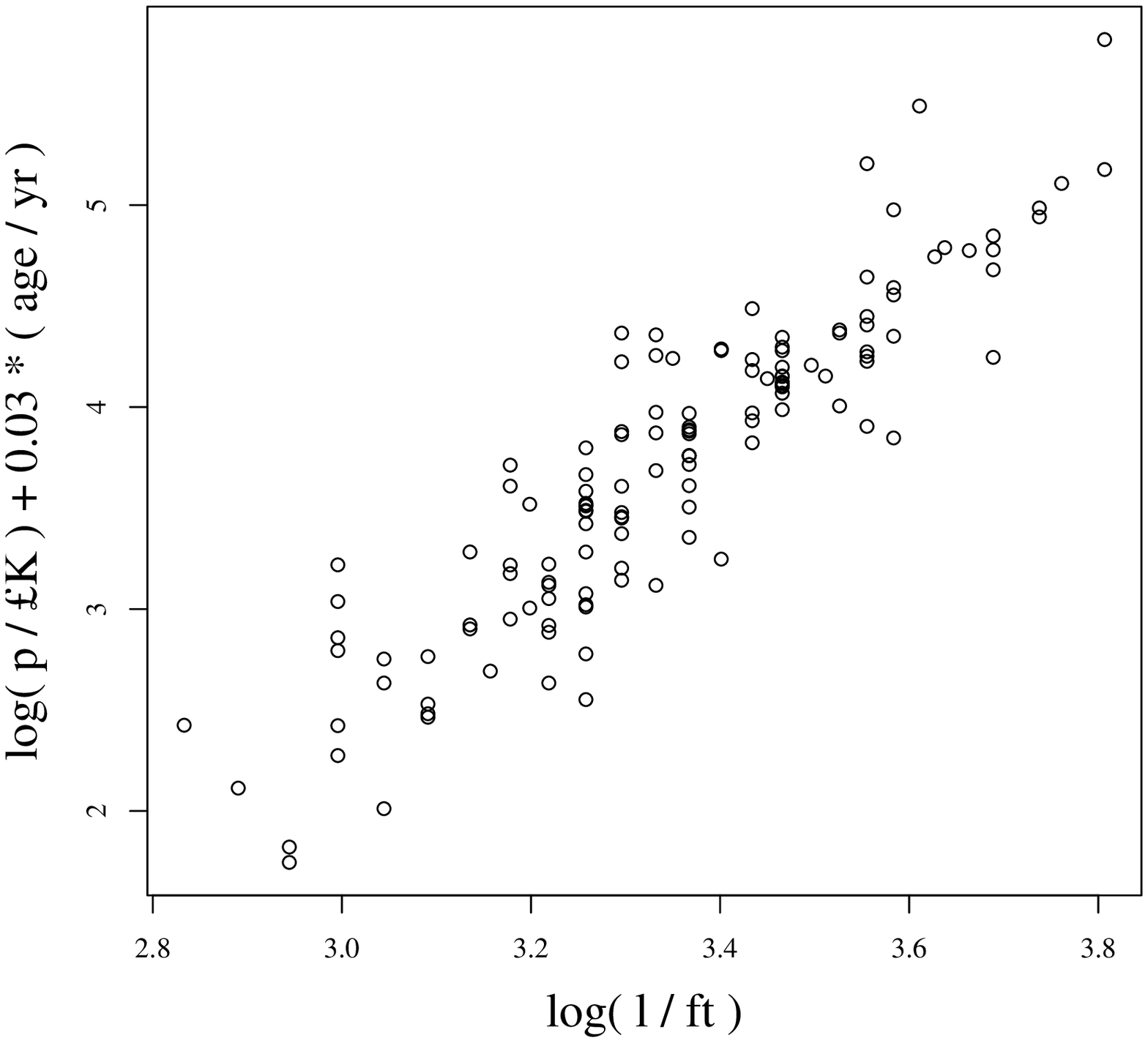}}

\centerline{Figure 4: how the price $p$ of secondhand yachts scales with their length $l$:}
 \centerline{(a) raw data $p \sim l^{3.8\pm0.2}, \;\sum R^2=0.71$; (b) age-adjusted $p \sim l^{3.5\pm0.1}e^{-0.03({\rm age / yr})}, \;\sum R^2=0.87$.}

\noindent The simplest way to look for a power law $y= C x^\beta$ is to take its logarithm,\footnote{We use natural logarithms, but any base will do.} so that $\log y =\beta\log x + \log C$. Rescalings of $x$ and $y$ are translations of the origin. A power law is seen as a straight line, and the scale-invariance property is that a straight line remains straight under changes of origin. Fig.\ 4(a) shows the prices of some used sailing yachts \cite{PBO}, and (because they were secondhand) Fig.\ 4(b) removes the best-fitting exponential decay of price with age, about 3\% of value per year. Both the high resulting power (price $p\sim l^{3.5}$) and the quality of the fit (nearly 90\% of variance) are striking. But notice that, for a given length and age, there's still typically an $e$-fold variation in prices---due not only to condition, design and desirability, but probably also to the enormous cost of the kit (rig, engine, electronics {\em etc.}) on board.

\vspace{0.1in}\noindent{\bf III.} To an applied mathematician, the Body Mass Index (BMI), your mass $m$  divided by the square of your height $h$, is a peculiar construct. It's clearly not invariant under self-similar scaling---why not the cube of height?---but then perhaps neither is a person's health: a six-foot man with the body shape of a toddler is probably not in good shape. In fact the reason for the square is empirical: among adults, levels of fat are maximally correlated with $mh^{-\beta}$ when $\beta\simeq2$ \cite{Quetelet}. But for children, where there's some genuine scaling going on, it seems highly inappropriate to use BMI, whose meaning will then be size-dependent. So, if $m\sim h^\beta$, what is $\beta$ for children?

An old argument for quadrupeds regards the animal as a beam supported by its legs, and derives the scaling law which places an upper bound on the proportionate sag of the beam \cite{Rash}. If we apply this instead to a human, regarding the pelvis as a supporting beam for the torso, one finds $7/3<\beta<8/3$, probably at the upper end of the range \cite{MacK}. For  boys aged 5-18 this matches the age-binned data rather well (Figure 5).\footnote{Growing girls have a power closer to three, due to their optimization for reproduction.}

\centerline{\includegraphics[keepaspectratio=true,
width=200pt]{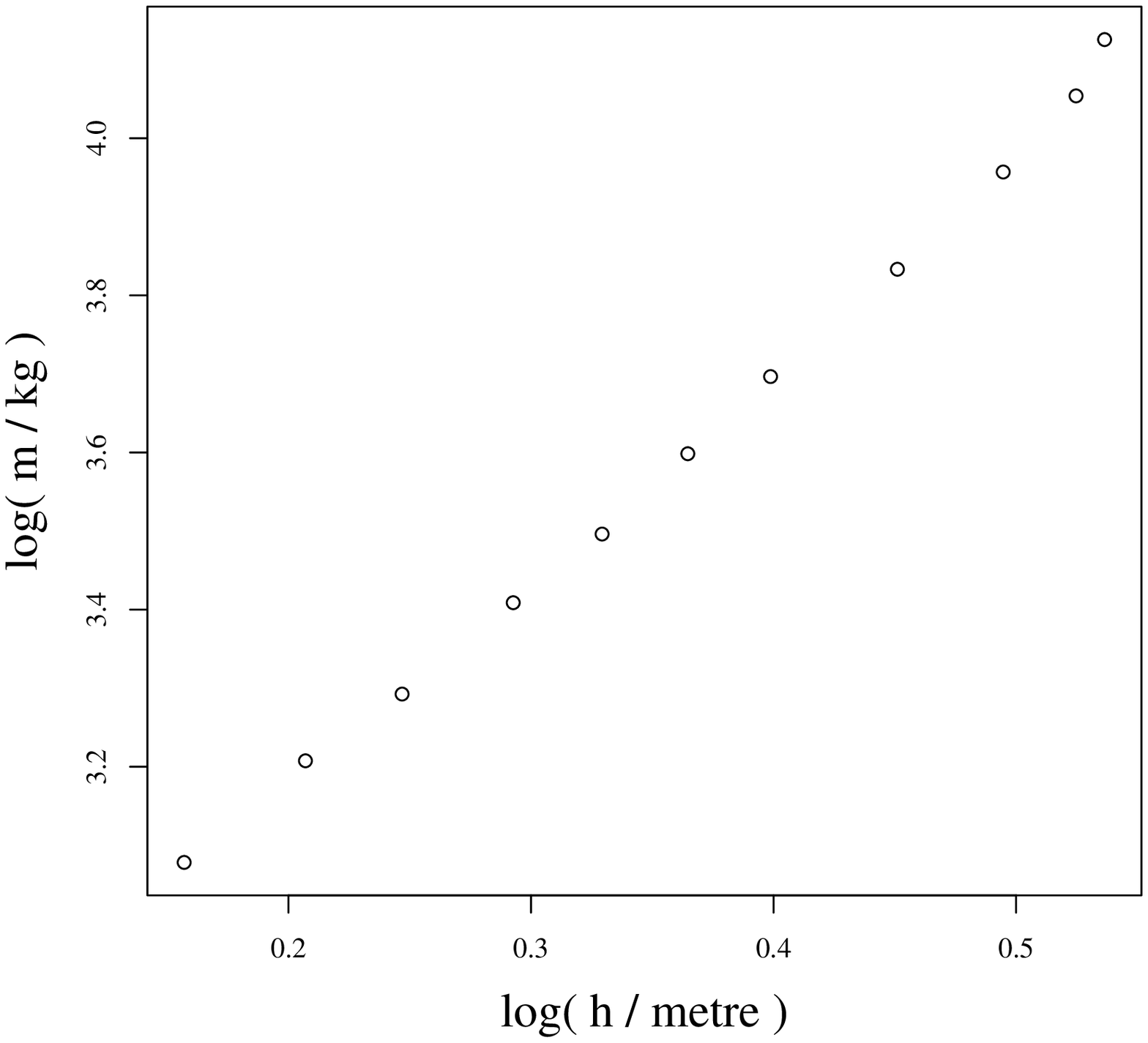}\hspace{0.6in}\includegraphics[keepaspectratio=true,
width=200pt]{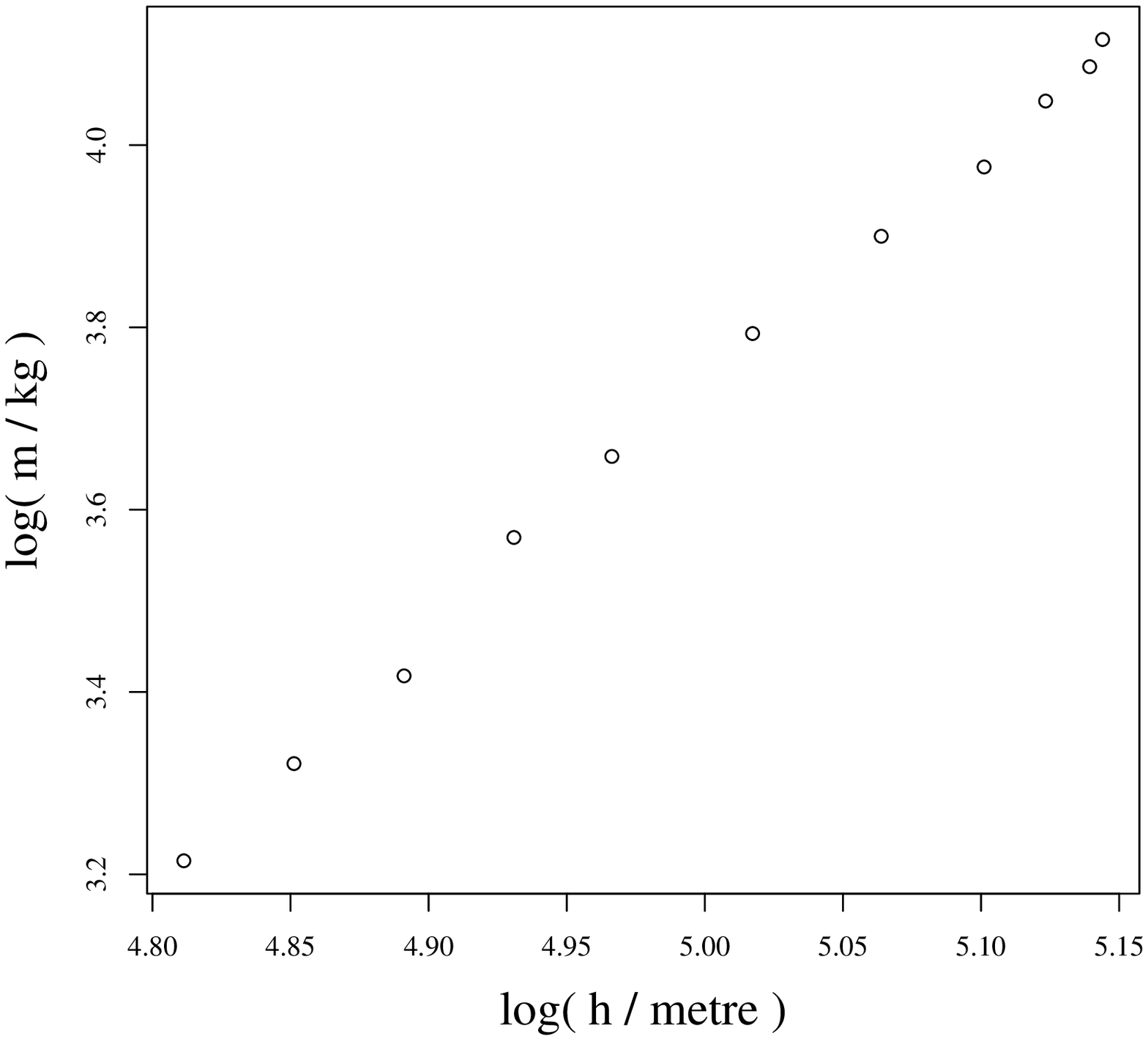}}

\centerline{Figure 5: how mass $m$ scales with height $h$ for boys in }
\centerline{(a) UK: $m \sim h^{2.70\pm0.05}$,  (b) Hong Kong: $m \sim h^{2.66\pm0.05}$; both $\sum R^2>0.995$.}

\vspace{0.1in}\noindent{\bf IV.} It's remarkable how well the basic structure of mammals scales: the dormouse and the bear in Fig.\ 6  share the same basic architecture, even though they are four orders of magnitude apart in mass. So what scaling laws might apply to mammals generally?
\vskip 0.2in
\centerline{
\hspace*{0in}\includegraphics[keepaspectratio=true,
width=135pt]{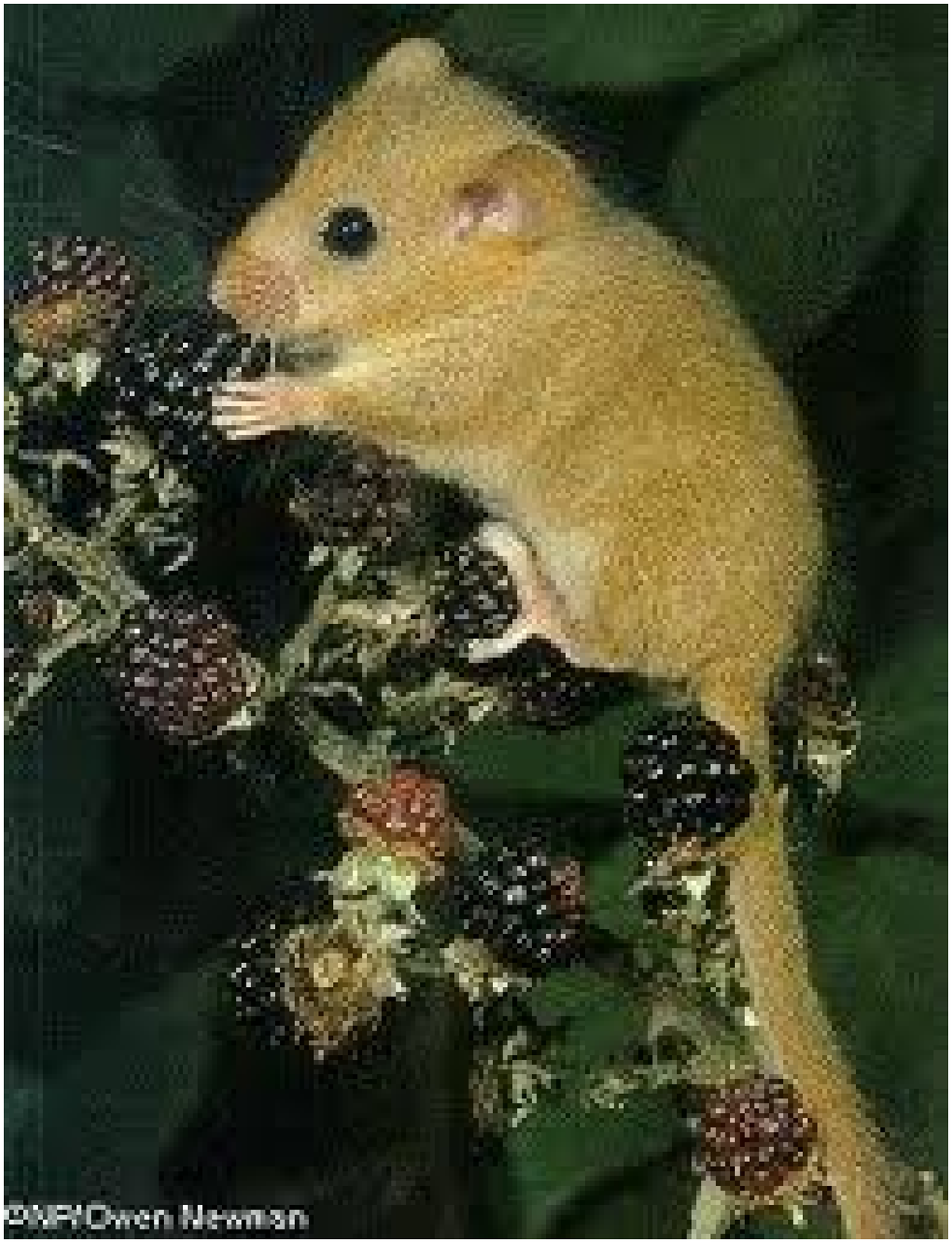}
\hspace*{0.5in}\includegraphics[keepaspectratio=true,
width=133pt]{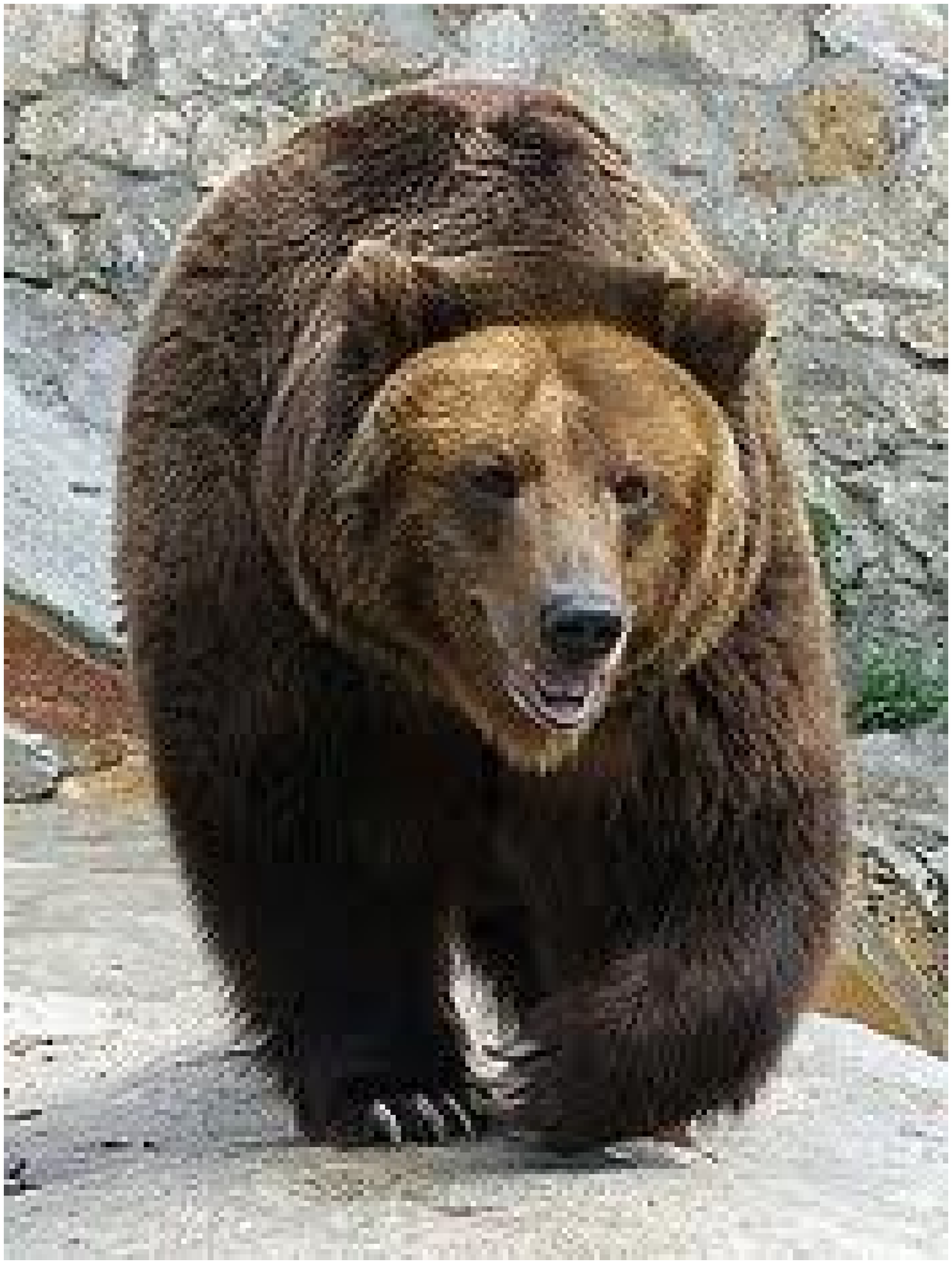}}

\centerline{Figure 6: a dormouse (mass $\sim20$g) and a bear ($200$kg).}

\noindent In his classic essay `On being the right size', J. B. S. Haldane noted that `You can drop a mouse down a thousand-yard mine shaft; and, on arriving at the bottom, it gets a slight shock and walks away, provided that the ground is fairly soft' \cite{Haldane}. For a first treatment of this, let's presuppose self-similar scaling,\footnote{Clearly the eyes \cite{eyes} and legs \cite{Haldane} scale differently from the torso. Biologists refer to non-self-similar scaling as `allometry', in contrast to `isometry'.} so that for an animal of length scale $l$ we have $m\sim l^3$. Air resistance is proportional to cross-sectional surface area (and thus to $l^2$) and to the square of speed $v$. Then, at terminal velocity, with air resistance balancing gravity, we have $l^2v^2\sim l^3$ and thus $v\sim l^{1/2}\sim m^{1/6}$. If the bear's terminal velocity is roughly 150mph, then that of the mouse is roughly 30mph. To do the analysis properly we should further consider the scaling of the process of sudden deceleration, but the effect is the same: as Haldane put it, `A rat is killed, a man is broken, a horse splashes'.

A classic scaling law for animals is that of metabolic rate $s$ (the animal's power consumption) with mass $m$. A simple self-similarity argument, with the animal losing heat in proportion to its surface area, gives $s\sim m^{2/3}$, but empirical data suggest $s\sim m^{3/4}$, `Kleiber's Law' \cite{Kleiber}. With modern data, of carefully measured {\em basal} metabolic rates, and for marsupials (which are less variable in their exploitation of new energy sources than other mammals), the fit is amazing (Fig.\ 7(a)) \cite{data}. One persuasive explanation \cite{WBE} derives Kleiber's law from the scaling of branching networks, which nature uses widely: think of alveoli, vascular networks, tracheoles, trees. Lungs, for example, are effectively neither three- nor wholly two-dimensional, rather filling the lung cavities with as much oxygen-grabbing potential as a branching network can manage. But (returning for a moment to flippancy) we could instead create a connection with Sect.\ III by  noting that if mammals' scaling is not self-similar but rather $m\sim l^{8/3}$, but their surface area is still dominated by terms proportional to $l^2$, then the heat-loss-through-surface argument immediately gives $s\sim l^2\sim m^{3/4}$ ---Kleiber's Law!

\vskip 0.2in
\centerline{\includegraphics[keepaspectratio=true,
width=200pt]{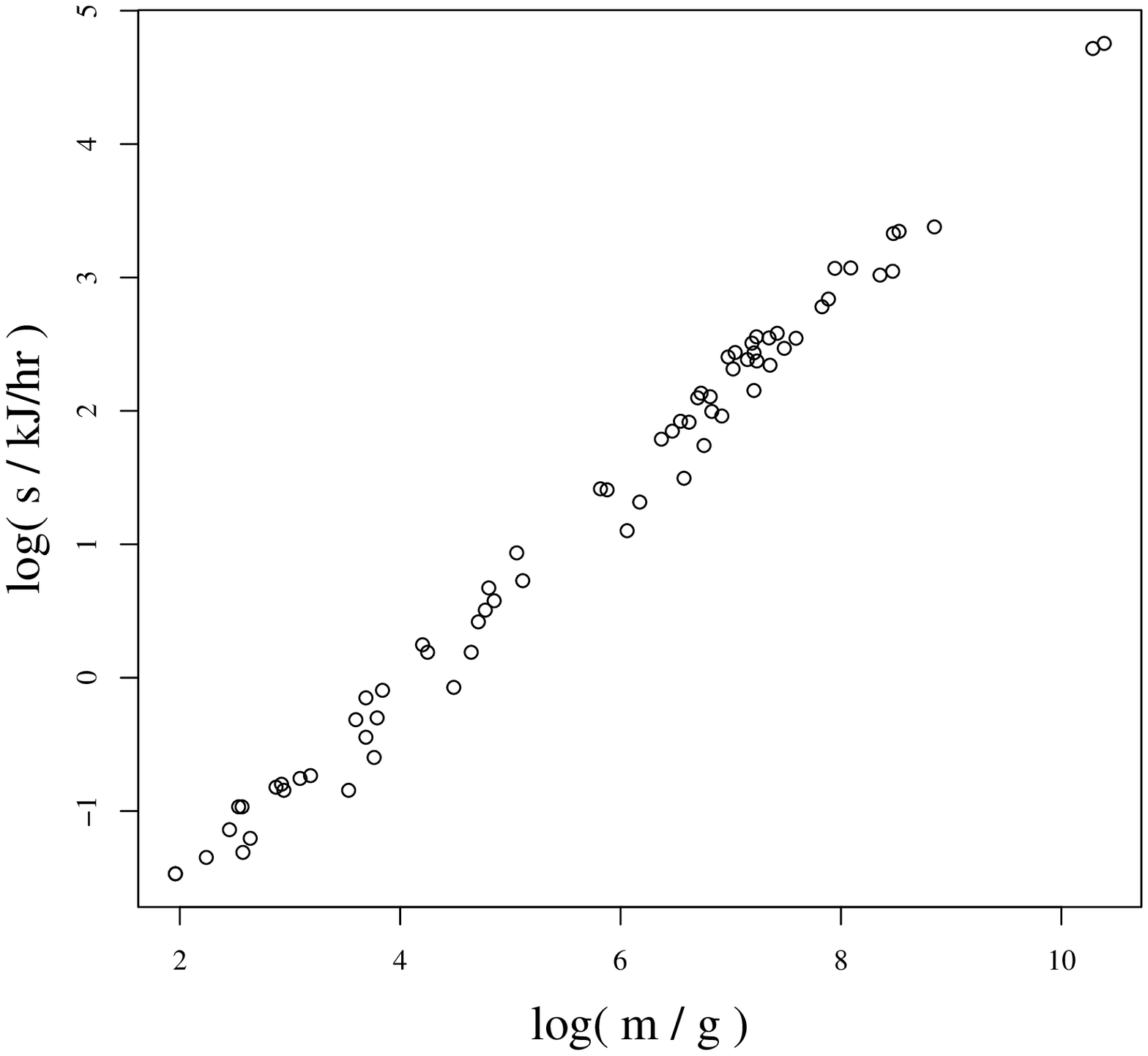}\hspace{0.6in}\includegraphics[keepaspectratio=true,
width=200pt]{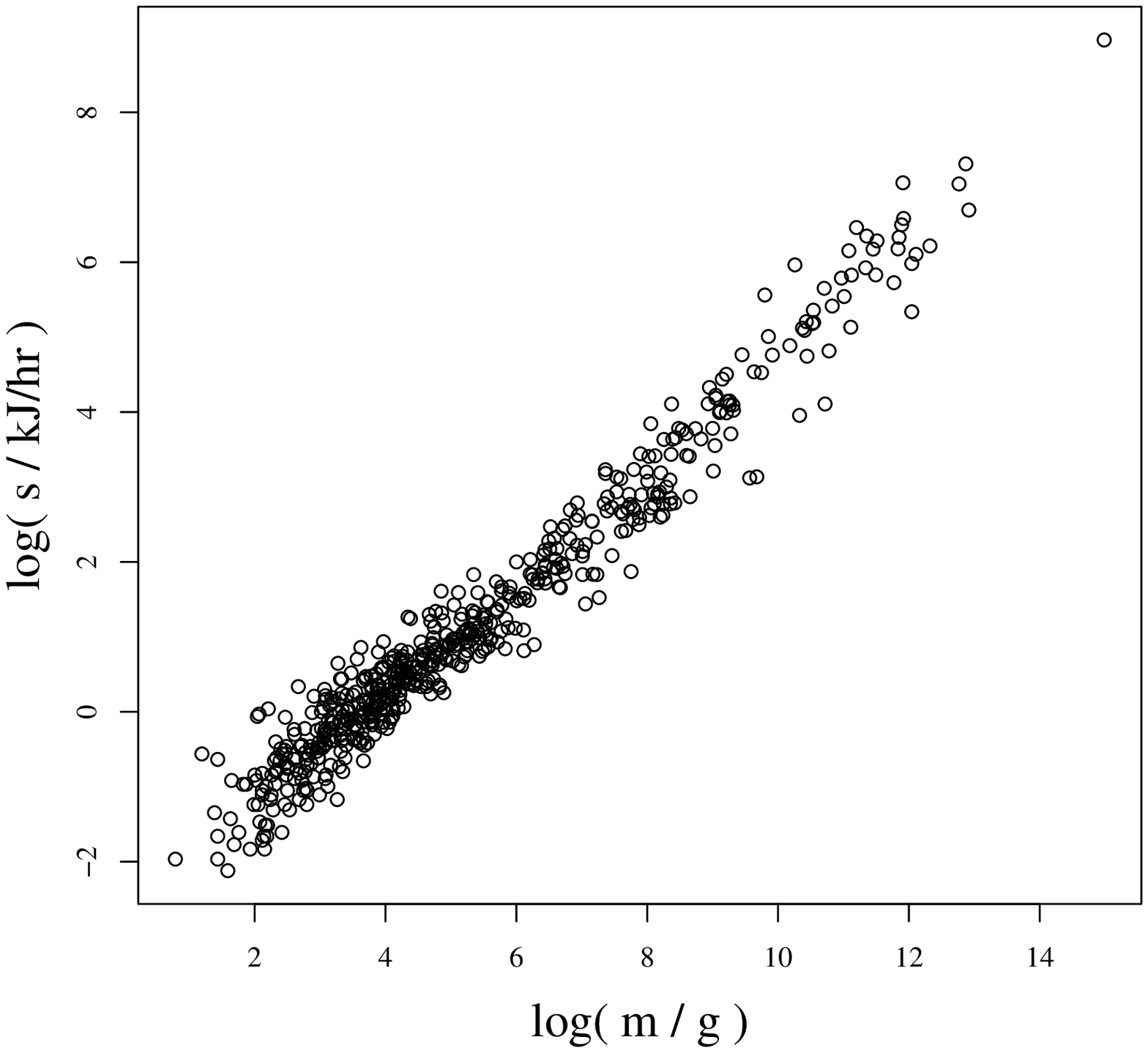}}

\centerline{Figure 7: how metabolic rate $s$ scales with mass $m$ for}\centerline{(a) Marsupials: $s \sim m^{0.75\pm0.01},\; \sum R^2=0.99$; (b) Eutheria: $s \sim m^{0.72\pm0.01},\; \sum R^2=0.96$.}

This is a good point at which to answer a fundamental question about data fitting: are we correct to fit a line to the logarithmic data, rather than a power to the untransformed, natural data? The question is really about how we weight the residuals---in the extreme cases, about how we treat the outliers. Suppose there were a species of mouse with ten times the metabolic rate expected for its size, and a species of bear with a rate 10\% greater than expected. Which is the outlier? If we believe that there is a scaling relation, a similarity in the architecture of creatures of such disparate size---and clearly we must do so if we're fitting a power---then it's the mouse that's the outlier, the stranger animal. This is treated correctly by the logarithmic data, where the mouse is about 24 times further from the fitted line than the bear. But in the natural data it's the bear that's the outlier, about 16 times further from the curve.\footnote{Nor does using multiplicative errors solve the problem, for (at a similar mass) we wish to give the same significance to a creature with double the expected BMR as to one with half the expected BMR. The logarithmic data do this, for $|\log 2|=|\log{1\over 2}|$, but merely ussing a multiplicative error weighting does not, for $2-1>1-{1\over 2}$.}

For the eutheria---the `good beasts', the non-marsupial mammals---the picture is not quite so clear (Fig.\ 7(b)), and indeed the subject of metabolic scaling has become surprisingly controversial in recent years \cite{BMR}. One group \cite{Kolo} attempted to explain the slight observable curvature in Fig.\ 7(b) by altering the model from $s\sim m^\beta$ and thus
$$ \log s = \alpha + \beta \log m$$
(where $\alpha$ and $\beta$ are the constants to be fitted by linear regression) to
$$ \log s = \alpha + \beta \log m + \gamma (\log m)^2\,.$$
This illustrates perfectly why all scientists need to be aware of dimensional issues, for  it is impossible to take the logarithm of a mass. `But,'  might come the response, `of course I can. This mouse weighs 20g. Look: I type 20 into my calculator, press ``log", and there you are'. But a mass has no intrinsic number, until we specify a unit. Here this is one gram, and what has actually been done is to take the logarithm of the dimensionless ratio of a mass $m$ to a mass unit $m_0=1g$.\footnote{This explains the precise form of the axis labels in Figs 4, 5, 7 \& 8. In Fig.\ 4, for example, `$\log( l/$ft$)$' means the logarithm of the ratio of the length $l$ to one foot.} Thus the model is actually
$$ \log \left({s\over s_0}\right) = \alpha + \beta \log \left({ m\over m_0}\right) + \gamma \log^2 \left({ m\over m_0}\right)\,,$$
and is no longer invariant under change of unit. If we replace $m_0$ by $m_0'$, then both $\alpha$ and $\beta$ are altered: with $\mu:=\log(m_0/m_0')$, we have $\alpha\mapsto\alpha+\beta\mu+\gamma\mu^2$ and $\beta\mapsto \beta+2\gamma\mu$. Scale invariance of the functional form has been lost, and in particular the linear term, the power $\beta$, is now a function of the choice of unit, and can always be set to zero by a correct choice of $m_0$, for (unlike a line) a parabola {\em does} have a preferred origin. The coefficient $\gamma$ of the quadratic term remains invariant, and is proportional to the change in gradient of the logarithmic data---but the mere fact of being able to fit a non-trivial quadratic is meaningless, for  a line segment is approximated arbitrarily well by a quadratic, provided its turning point is sufficiently far away (Figure 8) \cite{MacK2}.

\centerline{\includegraphics[keepaspectratio=true,
width=200pt]{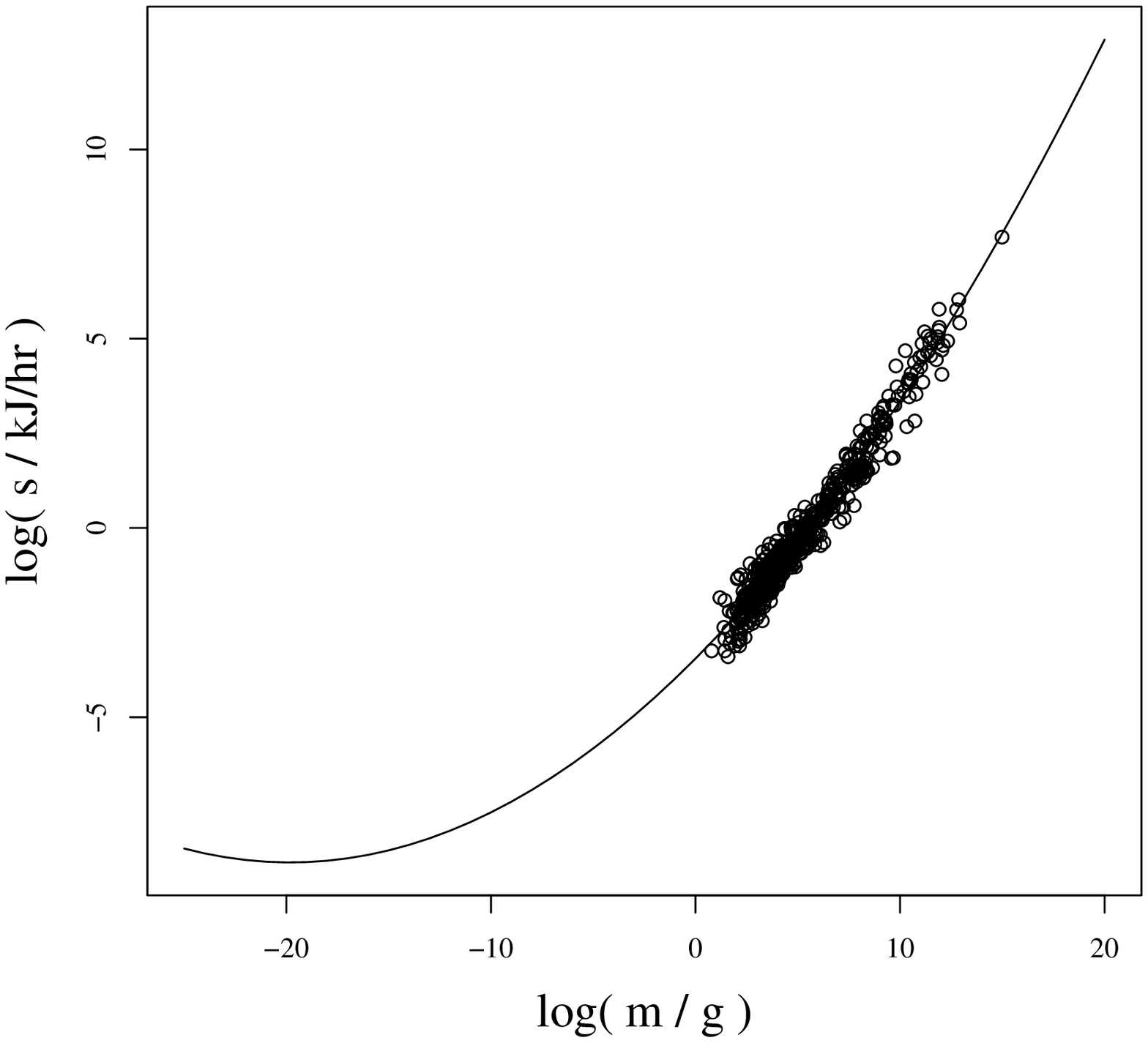}}

\centerline{Figure 8: the quadratic fitted to eutherian data.}

\vspace{0.1in}\noindent{\bf V.} The central role of the humble power among the functions used in applied mathematics is sometimes masked by the more prominent contributions of the exponential, trigonometric and even many of the special functions, the mathematics of which is more substantial. But all science---all applied maths---has units and dimensions, and so has dimensionless numbers, and so has scaling and thus power laws.

The initial ideas we have rambled through here lead on to all sorts of profound topics. Power law probability distribution functions and the processes which result in them make a lovely story in themselves \cite{pdfs}. One of the most mathematically-sophisticated subjects of physics, the theory of critical phenomena, has as its output universal classes of such power laws \cite{stanley}. Most recently, the idea that natural processes can enforce criticality (rather than needing to be fine-tuned) has been revolutionary \cite{bak}. We wish the reader joy in exploring these topics!

\end{document}